\documentclass[a4paper,english,fontsize=10pt,parskip=half,abstracton]{scrartcl}
\pdfoutput=1
\usepackage{babel}
\usepackage[utf8]{inputenc}
\usepackage[T1]{fontenc}
\usepackage[a4paper,left=20mm,right=20mm,top=30mm,bottom=30mm]{geometry}
\usepackage{amsmath}
\usepackage{amsthm}
\usepackage{amssymb}
\usepackage{enumerate}
\usepackage{aliascnt}
\usepackage[bookmarks=false,
            pdftitle={A counterexample to a conjecture of Kiyota, Murai and Wada},
            pdfauthor={Benjamin Sambale},
            pdfkeywords={Cartan matrices, eigenvalues, rationality},
            pdfstartview={FitH}]{hyperref}

\newaliascnt{Lem}{Thm}

\aliascntresetthe{Lem}
\newaliascnt{Prop}{Thm}

\aliascntresetthe{Prop}
\newaliascnt{Cor}{Thm}

\aliascntresetthe{Cor}
\newaliascnt{Conj}{Thm}
\newtheorem{Conj}[Conj]{Conjecture}
\aliascntresetthe{Conj}
\newaliascnt{Quest}{Thm}
\newtheorem{Quest}[Quest]{Question}
\aliascntresetthe{Quest}
\theoremstyle{definition}
\newaliascnt{Def}{Thm}

\aliascntresetthe{Def}
\newaliascnt{Ex}{Thm}

\aliascntresetthe{Ex}

\allowdisplaybreaks[1]

\renewcommand{\phi}{\varphi}

\newcommand{\N}{\operatorname{N}}

\newcommand{\ZZ}{\mathbb{Z}}

\newcommand{\FF}{\mathbb{F}}

\newcommand{\SL}{\operatorname{SL}}
\newcommand{\PSU}{\operatorname{PSU}}

\newcommand{\PSp}{\operatorname{PSp}}

\newcommand{\tr}{\operatorname{tr}}

\mathchardef\ordinarycolon\mathcode`\:  
 \mathcode`\:=\string"8000
 \begingroup \catcode`\:=\active
   \gdef:{\mathrel{\mathop\ordinarycolon}}
 \endgroup

\title{A counterexample to a conjecture\\ of Kiyota, Murai and Wada}
\author{Benjamin Sambale\footnote{Fachbereich Mathematik, TU Kaiserslautern, 67653 Kaiserslautern, Germany, 
\href{mailto:sambale@mathematik.uni-kl.de}{sambale@mathematik.uni-kl.de}}}
\date{\today}

\begin{document}
\frenchspacing
\maketitle
\begin{abstract}\noindent
Kiyota, Murai and Wada conjectured in 2002 that the largest eigenvalue of the Cartan matrix $C$ of a block of a finite group is rational if and only if all eigenvalues of $C$ are rational. We provide a counterexample to this conjecture and discuss related questions.
\end{abstract}


Let $B$ be a block of a finite group $G$ with respect to an algebraically closed field of characteristic $p>0$.
It is well-known that the Cartan matrix $C\in\ZZ^{l\times l}$ of $B$ is symmetric, positive definite, non-negative and indecomposable (here $l=l(B)$ is the number of simple modules of $B$). Let $ED$ (respectively $EV$) be the multiset of elementary divisors (respectively eigenvalues) of $C$. Note that these multisets do not depend on the order of the simple modules of $B$. Let $D$ be a defect group of $B$. Then the elementary divisors of $C$ divide $|D|$ and $|D|$ occurs just once in $ED$. 
On the other hand, the eigenvalues of $C$ are real, positive algebraic integers. By Perron-Frobenius theory, the largest eigenvalue $\rho(C)$ (i.\,e. the \emph{spectral radius}) of $C$ occurs with multiplicity $1$ in $EV$. Moreover, $\prod_{\lambda\in ED}\lambda=\det(C)=\prod_{\lambda\in EV}\lambda$. Apart from these facts, there seems little correlation between $EV$ and $ED$.

According to (the weak) Donovan's Conjecture, there should be an upper bound on $\rho(C)$ in terms of $|D|$. However, it can happen that $\rho(C)>l(B)|D|$. For example, if $B$ is the principal $2$-block of $G=\PSp(4,4).4$, a computation with GAP~\cite{GAP48} shows that $\rho(C)>7201>5\cdot 2^{10}=l(B)|D|$. 
This is even more striking than the observation $\tr(C)>l(B)|D|$ made in \cite{NS} for the same block. Conversely, $|D|$ cannot be bounded in terms of $\rho(C)$: for $p\ge 5$ the principal $p$-block of $\SL(2,p)$ satisfies $\rho(C)<4<p=|D|$ (see \cite[Example on p. 3843]{Wada8}). 

If $\lambda\in EV\cap\ZZ$, then $|D|/\lambda$ is an eigenvalue of $|D|C^{-1}\in\ZZ^{l\times l}$ and therefore it is an algebraic integer. This shows that $\lambda$ divides $|D|$. By a similar argument, $\lambda$ is divisible by the smallest elementary divisor of $C$. 
In \cite[Questions~1 and 2]{Wada5}, Kiyota, Murai and Wada proposed the following conjecture on the rationality of eigenvalues (see also \cite[Conjecture]{WadaEV}). 

\begin{Conj}[Kiyota-Murai-Wada]\label{wada}
The following assertions are equivalent:
\begin{enumerate}[(1)]
\item $EV=ED$.
\item $\rho(C)=|D|$.
\item $\rho(C)\in\ZZ$.
\item $EV\subseteq\ZZ$.
\end{enumerate}
\end{Conj}

Clearly, (1) $\Rightarrow$ (2) $\Rightarrow$ (3) $\Leftarrow$ (4) $\Leftarrow$ (1) holds and it remains to prove (3) $\Rightarrow$ (1). This has been done for blocks of finite or tame representation type (see \cite[Propositions~3 and 4]{Wada5}). For $p$-solvable $G$ we have (1) $\Leftrightarrow$ (2) $\Leftrightarrow$ (4) and $\rho(C)\le|D|$ (see \cite[Theorem~1]{Wada5}, \cite[Corollary~3.6]{Wada1} and \cite[Corollary 3.6]{Wada8}).
Other special cases were considered in \cite{KoshitaniYoshii2010,Wada2,Wada4,Yutaka}.
If $D\unlhd G$, then (1)--(4) are satisfied (see \cite[Proposition~2]{Wada5}). This holds in particular for the Brauer correspondent $b$ of $B$ in the normalizer $\N_G(D)$. 
In view of Broué's Abelian Defect Group Conjecture, Kiyota, Murai and Wada~\cite[Question~3]{Wada5} raised the following question.

\begin{Quest}[Kiyota-Murai-Wada]\label{Q}
If $D$ is abelian and $\rho(C)=|D|$, are $B$ and $b$ Morita equivalent?
\end{Quest}

It was proved in \cite{Wada2,KoshitaniYoshii2010} that the answer to \autoref{Q} is yes for principal $p$-blocks whenever $p\in\{2,3\}$.
However, the following counterexample shows not only that \autoref{wada} is false, but also that \autoref{Q} has a negative answer (for principal blocks) in general:

Let $B$ be the principal $5$-block of $G=\PSU(3,4)$. The Atlas of Brauer characters~\cite{BrauerAtlas} (or \cite{MOC}) gives
\[C=\begin{pmatrix}
10&10&5\\
10&13&6\\
5&6&7
\end{pmatrix}.\]
It follows that $EV=\bigl\{\frac{1}{2}(5+\sqrt{5}),\,\frac{1}{2}(5-\sqrt{5}),\,25\bigr\}$ and $ED=\{1,\,5,\,25\}$. Therefore, $\rho(C)=25=|D|$, but $EV\ne ED$. Moreover, $D$ is abelian since $|D|=25$, but $B$ cannot be Morita equivalent to $b$, since the eigenvalues of the Cartan matrix of $b$ are rational integers as explained above.

We do not know whether the implications (3) $\Rightarrow$ (2), (4) $\Rightarrow$ (1) or (4) $\Rightarrow$ (2) in \autoref{wada} might hold in general. Wada~\cite[Decomposition Conjecture]{Wada4} strengthened all three implications as follows.

\begin{Conj}[Wada]\label{w2}
There exist partitions $EV=E_1\sqcup\ldots\sqcup E_n$ and $ED=F_1\sqcup\ldots\sqcup F_n$ of multisets such that
\begin{itemize}
\item $|E_i|=|F_i|$ for $i=1,\ldots,n$.
\item $\prod_{\lambda\in E_i}\lambda=\prod_{\lambda\in F_i}\lambda$ for $i=1,\ldots,n$.
\item $\prod_{\lambda\in E_i}(X-\lambda)\in\ZZ[X]$ is irreducible for $i=1,\ldots,n$.
\item $\rho(C)\in E_1$, $|D|\in F_1$.
\end{itemize}
\end{Conj}

Again we found a counterexample: The group $\PSU(3,3)$ has a faithful $7$-dimensional representation over $\FF_3$. Let $G=\FF_3^7\rtimes\PSU(3,3)$ be the corresponding semidirect product, and let $B$ be the principal $3$-block of $G$. 
This group and its character table can be accessed as $\texttt{PrimitiveGroup}(3^7,35)$ and $\texttt{CharacterTable("P49/G1/L1/V1/ext3")}$ in GAP. In this way we obtain $9\in EV$, but $9\notin ED$. Obviously, this contradicts \autoref{w2}.

\section*{Acknowledgment}
I thank Thomas Breuer for some explanations about character tables in GAP. Moreover, I am grateful to Gabriel Navarro for getting me interested in counterexamples.
This work is supported by the German Research Foundation (project SA 2864/1-1).

\end{document}